\documentclass [reqno]{amsart}
\usepackage {times,mathptm}
\usepackage {amsmath,amssymb}
\usepackage {color,graphicx}
\usepackage [latin1]{inputenc}

\newcommand {\RR}{{\mathbb R}}
\newcommand {\ZZ}{{\mathbb Z}}

\newcommand {\preprint}[2]{preprint \discretionary {#1/}{#2}{#1/#2}}
\newcommand {\calM}{{\mathcal M}}

\DeclareMathOperator {\im}{im} 
\DeclareMathOperator {\val}{val} \DeclareMathOperator {\dist}{dist}

\DeclareMathOperator {\divi}{div}

\DeclareSymbolFont {mysymbols}{OMS}{cmsy}{m}{n}
\DeclareMathSymbol {\calI}{\mathalpha}{mysymbols}{`I}

\DeclareSymbolFont {mylargesymbols}{OMX}{cmex}{m}{n}
\DeclareMathSymbol {\dunion}{\mathop}{mylargesymbols}{"60}

\newcommand {\diunion}{\,\mbox {\raisebox{0.25ex}{$\cdot$} \kern-1.83ex $\cup$}
  \,}

\renewenvironment {enumerate}%
  {\rule{1mm}{0mm}\begin {oldenumerate}%
    \parskip1ex plus0.5ex \itemsep 0mm \parindent 0mm}%
  {\end {oldenumerate}}

\renewenvironment {itemize}%
  {\rule{1mm}{0mm}\begin {olditemize}%
    \parskip1ex plus0.5ex \itemsep 0mm \parindent 0mm}%
  {\end {olditemize}}

\theoremstyle {plain}
\newtheorem {theorem}{Theorem}[section]
\newtheorem {proposition}[theorem]{Proposition}
\newtheorem {lemma}[theorem]{Lemma}
\newtheorem {corollary}[theorem]{Corollary}
\newtheorem {lemmaanddef}[theorem]{Lemma and Definition}
\theoremstyle {definition}
\newtheorem {definition}[theorem]{Definition}
\theoremstyle {remark}
\newtheorem {remark}[theorem]{Remark}
\newtheorem {example}[theorem]{Example}

\newtheorem {notation}[theorem]{Notation}

\begin{document}

\title [Intersecting Psi-classes on tropical $\mathcal{M}_{0,n}$]{Intersecting Psi-classes on tropical $\mathcal{M}_{0,n}$}
\author {Michael Kerber and Hannah Markwig}
\address {Michael Kerber, Fachbereich Mathematik, Technische Universit\"at
  Kaiserslautern, Postfach 3049, 67653 Kaiserslautern, Germany}
\email {mkerber@mathematik.uni-kl.de}
\address {Hannah Markwig, Institute for Mathematics and its
Applications (IMA), University of Minne\-sota, Lind Hall 400, 207
Church Street SE, Minneapolis, MN 55455, USA}
\email{markwig@ima.umn.edu}
\thanks{The second author would like to thank the Institute for Mathematics and its Applications (IMA) in Minneapolis for hospitality and the German academic exchange service (DAAD) for financial support}

\begin{abstract}
We apply the tropical intersection theory as suggested by G.~ Mikhalkin and developed in detail by L.~Allermann
and J.~Rau to compute intersection products of tropical Psi-classes on the
moduli space of rational tropical curves. We show that in the case
of zero-dimensional (stable) intersections, the resulting numbers
agree with the intersection numbers of Psi-classes on the moduli
space of $n$-marked rational curves computed in algebraic
geometry.
\end{abstract}

\maketitle

  % Tropical Psi-Classes: Introduction

\section{Introduction}

A rational $n$-marked tropical curve is a metric tree with $n$
labeled leaves and without $2$-valent vertices. Those curves are parametrized by the
combinatorial structure of the underlying (non-metric) tree
and the length of each interior edge. The tropical moduli space
$\calM_{0,n}$ (the space which parametrizes these curves) has the
structure of a polyhedral complex, obtained by gluing several copies
of the positive orthant $\RR_{\geq 0}^{n-3}$
--- one copy for each $3$-valent combinatorial graph with
$n$ leaves (see section $2$).

Recently, G.~Mikhalkin (see \cite{M1}) introduced tropical Psi-classes on the moduli space of rational tropical curves:
for $k\in [n]$, the tropical Psi-class $\Psi_k$ is
the subcomplex of cones of $\calM_{0,n}$ corresponding to tropical curves
which have the property that the leaf labeled with the number $k$ is
adjacent to a vertex of valence at least $4$ (see definition
\ref{def-Psi}).

The aim of this article is to apply the concepts of tropical
intersection theory suggested by G.~Mikhalkin and developed in detail by L.~Allermann and J.~Rau (\cite{M2}, \cite{AR})
to compute the intersection products of an arbitrary number of these
Psi-classes.

In order to do this, we first recall the embedding of the moduli
space $\calM_{0,n}$ of $n$-marked rational tropical curves into some
real vector space $Q_n$ (and other preliminaries) in section $2$. On this space $Q_n$, we construct in
section $3$ a tropical rational function $f_k$ for
all $k\in [n]$ with the property that the Cartier divisor of the
restriction of $f_k$ to (the embedding of) $\calM_{0,n}$ is (a
multiple of) the $k$-th Psi-class $\Psi_k$. We use this description
in section $4$ to compute the weights on the maximal cones of the
tropical fan obtained by intersecting an arbitrary number of
tropical Psi-classes. As a special case, we compute the weights of
($0$-dimensional) intersections of $n-3$ tropical Psi-classes
--- they agree with the $0$-dimensional intersection product of
$n-3$ Psi-classes on the moduli space of rational
$n$-marked curves computed in algebraic geometry.\\

The authors would like to thank L.~Allermann, I.~Ciocan-Fontanine,
A.~Gathmann and J.~Rau for useful conversations.
We also thank an anonymous referee for suggesting a more elegant proof for the identity in lemma \ref{applem1}.

  % Tropical Psi-Classes: Erster Teil - oder: Wo die Notationen eingefuehrt werden...

%{\bf ToDo:}
%\begin{itemize}
%   \item Prove Lemma \ref{lem-positivbasis}
%   \item Prove Proposition \ref{prop-basiccase}
%   \item It might be useful to use set-notation $e_S$ instead of
%   double index notation $e_{ij}$.
%\end{itemize}

\section{Preliminaries}

In the sequel, $n$ will always denote an integer greater than $2$.

An $n$-marked (rational) abstract tropical curve is a metric tree
$\Gamma$ (that is, a tree together with a length function assigning
to each non-leaf edge a positive real number) without
$2$-valent vertices and with $n$ leaves, labeled by numbers
$\{1,\ldots ,n\}$ (see \cite{GaMa}, definition 2.2). The space
$\calM_{0,n}$ of all $n$-marked tropical curves has the structure of
a polyhedral fan of dimension $n-3$ obtained from gluing copies of
the space $\RR_{>0}^k$ for $0\leq k\leq n-3$ --- one copy for each
combinatorial type of a tree with $n$ leaves and exactly $k$ bounded
edges. Its face lattice is given by $\tau\prec\sigma$ if and only if
the tree corresponding to $\tau$ is obtained from the tree
corresponding to $\sigma$ by contracting bounded edges. For more
details, see \cite{BHV}, section 2, or \cite{GaMa}, section 2.

% has the structure of a fan of dimension $n-3$
%whose $k$-dimensional cells correspond to
%
%
%(see GM/BHV).
%
%
%
%An \emph{abstract tropical curve} is a connected rational graph
%$\Gamma$ (that is, a tree) whose vertices have valence at least $3$.
%Unbounded edges (also called ends or leaves) are allowed. The
%bounded edges are equipped with a positive length. A $n$-marked
%abstract tropical curve is a tuple $C=(\Gamma,x_1,\ldots,x_n)$ where
%$\Gamma$ is an abstract tropical curve and $x_1,\ldots,x_n$ are
%distinct unbounded edges. (For more details, see \cite{GM}, definition 2.2.)
%The space $ \calM_{0,n} $ is defined to be the \emph{space of
%all $n$-marked abstract tropical curves} (modulo isomorphism) with
%exactly $n$ leaves.
%The \emph{(combinatorial) type} of a marked abstract tropical curve
%$(\Gamma,x_1,\ldots,x_n)$ is the data left when dropping the
%information about the lengths of the bounded edges. Lemma 2.10 of
%\cite{GM} says that there are only finitely many combinatorial
%types. The subset of curves of type $\alpha$ in $ \calM_{0,n}
%$ is the interior of a cone in $\RR^k$ (where $k$ is the number of
%bounded edges), given by the inequalities that all lengths are
%positive --- that is, it is the positive orthant of $\RR^k$. Example
%2.13 of \cite{GM} describes how these cones are glued locally in $
%\calM_{0,n} $. (This construction is equal to the \emph{space of (phylogenetic) trees} $\mathbb{T}_n$ in \cite{BHV}, page 9-13.)

In order to recall how $\calM_{0,n}$ can be embedded from \cite{GKM}, we need
the following notations:

Let $\mathcal{T}:=\{S\subset[n]\; :\; |S|=2\}$ denote the set of
two-element subsets of $[n]:=\{1,\cdots ,n\}$. Consider the space
$\RR^{\binom{n}{2}}$ indexed by the elements of $\mathcal{T}$ and
define a map to this space via

\begin{eqnarray*}
 \Phi_n:  \RR^n & \longrightarrow & \RR^{\binom{n}{2}}\\
          a & \longmapsto &
         (a_i+a_j)_{{\{i,j\}\in\mathcal{T}}}.
\end{eqnarray*}

Let $Q_n$ denote the quotient vector space
$\RR^{\binom{n}{2}}/\im(\Phi_n)$, which has dimension
$\binom{n}{2}-n$.

Furthermore, we define a map
\begin{eqnarray*}
\varphi_n: \calM_{0,n} & \longrightarrow & \RR^{\binom{n}{2}}\\
C & \longmapsto & \dist(\{i,j\})_{\{i,j\}\in\mathcal{T}}
\end{eqnarray*}
where $\dist(\{i,j\})$ denotes the sum of the lengths of all bounded
edges on the (unique) path between the leaf marked $i$ and the leaf
marked $j$.

%\begin{eqnarray*}
%  \varphi_n: \calM_{0,n} & \longrightarrow  & \RR^{\binom{n}{2}}\\
%             (\Gamma,x_1,\ldots,x_n)  & \longmapsto      & (\dist_\Gamma(x_i,x_j))_{\{i,j\}\in\mathcal{T}}
%\end{eqnarray*}

%where $\dist_\Gamma(x_i,x_j)$ denotes the distance between the
%\sout{unbounded edges (or }leaves $x_i$ and $x_j$, that is, the sum
%of the lengths of all edges on the (unique) path leading from $x_i$
%to $x_j$.

\begin{theorem}
 Using the map $\varphi_n$, $\calM_{0,n} $ can be
embedded as a tropical fan into $Q_n$. 
\end{theorem}
For a proof, see theorem 3.4 of \cite{GKM}
or theorem 3.4 of \cite{SS}.

Note that $\calM_{0,n} $ is a marked fan (see
definition 2.12 of \cite{GKM}): Let $\tau$ be a cone and $C$ be the
corresponding tropical curve where all lengths of bounded edges are
chosen to be one. Then $\tau$ is generated by the rays
$\varphi_n(C_i)$, where $C_i$ denotes a curve obtained from $C$ by
shrinking all but one bounded edge to length $0$. In particular, $\calM_{0,n}$ is a simplicial fan and the $\varphi_n(C_i)$ form a basis for the span of the cone. They even form a unimodular basis which follows from proposition 5.4 of \cite{GiMa} (note that there, a different lattice and a different embedding of $\calM_{0,n} $ is used).

\begin{notation}\label{def-vk}
For each subset $I\subset [n]$ of cardinality $1<|I|<n-1$,
define a vector $v_I\in\RR^{\binom{n}{2}}$ via
$$(v_{I})_{T}:=\begin{cases}
                 1, & \text{ if } |I\cap T|=1\\
                 0, & \text{ otherwise.}
              \end{cases}$$
Note that $v_I$ is the image under $\varphi_n$ of a tree with one
bounded edge of length one, the marked ends with labels in $I$ on
one side of the bounded edge and the marked ends with labels in
$[n]\setminus I$ on the other, hence $v_I=v_{[n]\setminus I}$. 

For
$k\in [n]$, we define $V_k:=\{v_I\; :\; k\not\in I \text{ and }
|I|=2\}$.
\end{notation}

\begin{lemma}\label{lem-generatingSet}
% For each subset $I\subset [n]$ with $1<|I|<n-1$, define a vector
% $v_I\in Q_n$ by
% $$(v_{I})_{T}:=\begin{cases}
%                   1, & \text{ if } |I\cap T|=1\\
%                   0, & \text{ otherwise}
%                \end{cases}$$
% For each $k\in [n]$, denote by $V_k$ the set $V_k:=\{v_I\; :\;
% k\not\in I \text{ and } |I|=2\}$. \\
% Then the linear span of $V_k$ equals $Q_n$.
% For any $k\in [n]$, the set $V_k$ (from definition
%\ref{def-vk}) spans
%$Q_n=\RR^{\binom{n}{2}}/\im(\Phi_n)$.
For any $k\in [n]$, the linear span of the set $V_k$ equals $Q_n=\RR^{\binom{n}{2}}/\im(\Phi_n)$.
\end{lemma}

\begin{proof}

 We prove that for $S=\{s_1,s_2\}\in\mathcal{T}$, the
 $S$-th standard unit vector $e_S\in \RR^{\binom{n}{2}}$ is the sum of a linear
combination of elements in $V_k$ and an element of $\im(\Phi_n)$.

First assume that $k\not\in S=\{s_1,s_2\}$. Then it follows
immediately from the definitions of $v_S$ and $\Phi_n$ that
$e_S=(-v_S+\Phi_n(e_{s_1}+e_{s_2}))/2$, where $e_{s_i}$
denotes the $s_i$-th unit vector in $\RR^n$.\\

Now assume that $S=\{s_1,k\}$. We claim that
\begin{equation}\label{equation}
e_S=\frac{1}{2}\left(\sum_{I\in \mathcal{T}\; :\; I\cap S=\{s_1\}}v_I -
\Phi_n(a)\right), \end{equation} where $a\in\RR^n$ is the vector with
entries
$$a_i:=\begin{cases}
         n-4, & \text{ if } i=s_1\\
         0, & \text{ if } i=k\\
         1 & \text{ otherwise }
       \end{cases}$$

Check this equality in each component $T=\{t_1,t_2\}$.\\
The entry there is equal to one if and only if $|I\cap T|=1$ --- note that $s_1\in I$.\\

If $S\cap T=\emptyset$, then $\Phi_n(a)_T=2$ and we have $|I\cap
T|=1$ iff $I$ contains one element of $T$.

If $S\cap T=\{k\}$, then $\Phi_n(a)_T=1$ and we have $|I\cap T|=1$
iff $I=\{s_1\}\cup T\setminus\{k\}$.

If $S\cap T=\{s_1\}$, then $\Phi_n(a)_T=n-3$ and we have $|I\cap
T|=1$ iff $I\not = T$.

If $S\cap T=\{s_1,k\}$, then $\Phi_n(a)_T=n-4$ and we have $|I\cap
T|=1$ for all $(n-2)$ choices of $I$.

It follows that for $T\not =S$, we have
$$\left(\sum_{I\in \mathcal{T}\; :\; I\cap S=\{s_1\}}v_I -
\Phi_n(a)\right)_T=0$$

and for $T=S$, we have
$$\left(\sum_{I\in \mathcal{T}\; :\; I\cap
S=\{s_1\}}v_I - \Phi_n(a)\right)_T=n-2-(n-4)=2,$$ hence equation
(\ref{equation}) holds.
\end{proof}

\begin{lemma}\label{lem-balancing}
The sum over all elements $v_S\in V_k$ (from notation \ref{def-vk}) is an element in
$\im(\Phi_n)$, hence
$$\sum_{v_S\in V_k}v_S=0\in Q_n.$$
\end{lemma}
\begin{proof}
If $k\not\in T$, then there are $\binom{n-3}{1}\binom{2}{1}$ two-element subsets $S\in [n]\setminus\{k\}$ with the property that
$S\cap T=1$, hence $(\sum v_s)_T=2(n-3)$ in this case. If $k\in T$,
then the number of subsets $S\in [n]\setminus\{k\}$ satisfying
$|S\cap T|=1$ equals $n-2$, hence $(\sum v_s)_T=n-2$ in this case.
It follows that
$$\sum_{v_S\in V_k} v_{S} =\Phi_n(n-3,\ldots, n-3,1,n-3,\ldots ,n-3)$$ with
entry $1$ at position $k$.

\end{proof}

\begin{lemmaanddef}\label{rem-basisdarstellung}
Every element $v\in Q_n$
has a unique representation $$v=\sum_{v_{S}\in V_k}
\lambda_{S}v_{S}$$ with $\lambda_{S}\geq 0$ for all $S$ and
$\lambda_S=0$ for at least one $S\in\mathcal{T}$. We will call such
a representation in the future a \emph{positive representation of
$v$ with respect to $V_k$}.
\end{lemmaanddef}

\begin{proof}
As $|V_k|={\binom{n-1}{2}}=\dim(Q_n)+1$ (see notation \ref{def-vk}), the vectors of $V_k$ subdivide $Q_n$ into a fan whose $\dim(Q_n)+1$ top-dimensional cones are spanned by a choice of $\dim(Q_n)$ vectors of $V_k$. Each $v$ lies in a unique cone and its positive represenation is given by the linear combination of the spanning vectors of the cone.
Given a representation
$$v=\sum_{v_{S}\in V_k} \lambda_{S}v_{S}$$ with $\lambda_{S}> 0$ for
all $S$, the unique positive representation with respect to $V_k$
can be found by subtracting $\sum_{v_{S}\in V_k}
(\min_S{\lambda_{S}})v_{S}$.
\end{proof}
\begin{remark}\label{rem-basisdarstellung2}
It follows that a map from $V_k$ to $\RR_{\geq 0}$
gives rise to a well-defined convex piecewise-linear function on the
space $Q_n$ via $f(\sum\lambda_{S}v_{S}):=\sum \lambda_{S}
f(v_{S})$.
\end{remark}

\begin{lemma}\label{basisrepresentation}
Let
%$I\;\diunion J=[n]$
$I\subset [n]$ with $1<|I|<n-1$ and assume without restriction
that $k\notin I$. Then a positive representation of $v_I \in Q_n$
with respect to $V_k$ (as in definition \ref{rem-basisdarstellung} and notation \ref{def-vk}) is
given by
$$v_{I}=\sum_{S\subset I, v_{S}\in V_k}v_S.$$
\end{lemma}
\begin{proof}
 Let $|I|=m$, $I=\{i_1,\ldots,i_m\}$. We claim
that
\begin{displaymath}
v_{I}=\left(\sum_{S\subset I, v_{S}\in V_k}v_S \right)- (m-2)\cdot \Phi_n(e_{i_1}+\ldots+e_{i_m}).
\end{displaymath}
Check this equality in each component $T=\{t_1,t_2\}$. If $T\subset
I$, then $(v_I)_T=0$. There are $m-2$ choices for $v_S$ such that
$S$ contains $t_1$ and not $t_2$, and the same number of
choices such that $S$ contains $t_2$ and not $t_1$. Hence the first
sum of the right hand side contributes $2(m-2)$. As
\begin{displaymath}\Phi_n(e_{i_1}+\ldots+e_{i_m})_T=2\end{displaymath}
we get $0$ altogether. If $|T\cap I|=1$, then $(v_I)_T=1$. On the
right hand side, there are $m-1$ choices of $S$ such that $S$
contains $T\cap I$, and
\begin{displaymath}\Phi_n(e_{i_1}+\ldots+e_{i_m})_T=1.\end{displaymath}
If $T\cap I=\emptyset$,  both sides are equal to $0$.
\end{proof}

\section{Psi-classes as divisors of rational functions}

Let us start by reviewing some of the tropical intersection theory from \cite{AR}.
A \emph{cycle $X$} is a balanced, weighted, pure-dimensional, rational and 
polyhedral fan in $\RR^n$. The integer weights assigned to each top-dimensional cone $\sigma$ 
are denoted by $\omega(\sigma)$. By $|X|$, we denote the union of all cones of $X$ in $\RR^n$.
\emph{Balanced} means that the weighted
sum of the primitive vectors of the facets $\sigma_i$ around a cone $\tau \in X$ of codimension $1$
$$
  \sum_{i} 
    \omega(\sigma_i) u_{\sigma_i / \tau}
$$
lies in the linear vector space
spanned by $\tau$, denoted by $V_\tau$.
Here, a \emph{primitive vector $u_{\sigma_i/\tau}$ of $\sigma_i$ modulo
$\tau$} is a integer vector in $\ZZ^n$ that points from $\tau$ towards
$\sigma$ and fulfills the primitive condition: The lattice
$\ZZ {u}_{\sigma_i / \tau} + (V_\tau \cap \ZZ^n)$ must be equal to the 
lattice $V_{\sigma_i} \cap \ZZ^n$. Slightly differently, in \cite{AR}
the class of $u_{\sigma_i / \tau}$ modulo $V_\tau$ is called primitive
vector and $u_{\sigma_i / \tau}$ is just a representative of it. \\
Cycles are only considered up to refinements, i.e.\ we will consider two cycles equivalent if they have a common refinement.

A \emph{(non-zero) rational function on
$X$} is a continuous piece-wise linear function $\varphi : |X| \rightarrow \RR$ 
that is linear with rational 
slope on each cone. 
The \emph{Weil-divisor of $\varphi$ on $X$}, denoted by
$\divi(\varphi)$, is the
balanced subcomplex (resp. subfan) of $X$ defined in construction 3.3. of \cite{AR}, 
namely
the codimension one skeleton of $X$ together with weights $\omega(\tau)$ for each cone
$\tau \in X$ of codimension $1$. 
These weights are given by the formula
$$\omega(\tau)=\sum_i \varphi\left(\omega(\sigma_i)
u_{\sigma_i/\tau}\right)- \varphi\left(\sum_i \omega(\sigma_i)
u_{\sigma_i/\tau}\right),$$
where the sum goes again over all top-dimensional neighbours of $\tau$.

Now let us repeat the definition of tropical Psi-class.

\begin{definition}[see \cite{M1}, definition 3.1.]\label{def-Psi} For
%Following  G.~Mikhalkin (\cite{M}, definition 3.1), we define for
$k\in [n]$, the tropical Psi-class $\Psi_k\subset\calM_{0,n}$ is defined to be the weighted fan consisting of those closed $(n-4)$-dimensional cones that correspond to
tropical curves with the property that the leaf marked with the
number $k$ is adjacent to a vertex with valence $4$. The
weight of each cone is defined to be equal to one.
% one $4$-valent vertex adjacent to the leaf marked $k$. The
%weight of each cone is one.
\end{definition}
A motivation for this definition is given in \cite{M1}. Another motivation is that if one evaluates the classical $\psi_i$ on 1-strata of $\overline{M}_{0,n}$, i.e.\ on rational curves whose dual graphs are 3-valent except for one 4-valent vertex; we get $0$ if the $4$-valent vertex is not adjacent to the leaf $i$. To the author's knowledge, E.~Katz is about to prepare a preprint that explains more about the connection of tropical and classical Psi-classes.

In a recent article, G.~Mikhalkin defines an embedding of the space
$\calM_{0,n}$ as a tropical fan into
%(\cite{M}, theorem 1) gives a
%different embedding of the space $\calM_{0,n}$ as a tropical fan in
 $\RR^{\binom{n}{2}\binom{n-2}{2}}$,
  with the property that the tropical Psi-class $\Psi_k$ has the structure of a
  tropical subfan, i.e.\ satisfies the balancing condition (see \cite{M1},
theorem 3.1 and proposition 3.2).

Here, we prefer to work with the embedding of $ \calM_{0,n}$ in $Q_n$ as described in section 2.
We will see later in this section that
$\varphi_n(\Psi_k)\subset Q_n$ is a tropical fan --- this will
follow directly from the proof of proposition \ref{prop-basiccase},
which states that (a multiple of) $\varphi_n(\Psi_k)$ is the
Weil-divisor associated to a rational function.

By abuse of notation we will in the following not distinguish between $\Psi_k$ and $\varphi_n(\Psi_k)\subset Q_n$.

% \begin{lemma}\label{def-fk}
% Let $k\in [n]$ be a fixed number. Then the function $f_k$ given as
% the extension of the map $V_k\ni v_S\mapsto 1$ to the piecewise
% linear convex real-valued function on $Q_n$ defined by
% $$f(v)=f(\sum_{v_S\in V_k} \lambda_Sv_S) \mapsto \sum_{v_S\in V_k}
%   {\lambda_S}f(v_S),$$
% (where $\sum_{v_S\in V_k} \lambda_Sv_S$ is a positive representation
% of $v$ (see remark \ref{rem-basisdarstellung})) is linear on each
% cone of $\calM_{0,n}$.
% \end{lemma}
\begin{notation}\label{def-fk}
For any $k\in[n]$, let $f_k$ be the extension of the map $V_k\ni v_{S} \mapsto 1$ to
$Q_n$ (see notation \ref{def-vk} and remark \ref{rem-basisdarstellung2}).
\end{notation}
\begin{lemma}\label{lem-fklinear}
The map $f_k$ is linear on each cone of $\calM_{0,n}$.
\end{lemma}
\begin{proof}
Let $\tau$ be a cone and $C$ the tropical curve of the combinatorial
type corresponding to $\tau$ with all lengths equal to one. The cone
$\tau$ is generated by vectors $v_I$ corresponding to curves where
all but one bounded edge of $C$ are shrunk to length $0$. Assume
$\tau$ is generated by $v_{I_1},\ldots,v_{I_r}$, and let $k\notin
I_i$ for all $i$. Then each point $p$ in $\tau$ is given by a linear
combination $p=\sum_{i=1}^r \mu_i v_{I_i}$, where the $\mu_i$ are
non-negative. We can find a positive representation
\begin{displaymath} v_{I_i}=\sum_{v_S\in V_k} \lambda_{i,S} v_S\end{displaymath} for each $v_{I_i}$
using lemma \ref{basisrepresentation} (where each $\lambda_{i,S}$ is
either $1$ or $0$, depending on whether $S\subset I_i$ or not). We
claim that \begin{equation}\label{eq-pos-rep}p= \sum_{v_S\in
V_k}\big(\sum_{i=1}^r \mu_i \lambda_{i,S}\big) v_S\end{equation} is
a positive representation of $p$ with respect to $V_k$. It is
obvious that the $\sum_{i=1}^r \mu_i \lambda_{i,S}$ are
non-negative. It remains to show that there is at least one $S$ such
that $\sum_{i=1}^r \mu_i \lambda_{i,S}=0$. 
Let $a,b\in [n]$ be leaves in different
connected components of $C\setminus k$ 
(where $C\setminus k$ denotes the graph produced from $C$ by removing the closure of the unbounded edge labeled $k$, i.e.\ including the end vertex of $k$.)
There are at least two such connected components,
due to the fact that $k$ is adjacent to an at least $3$-valent vertex.
% 
% 
% 
% 
% The marked end $k$ is
% adjacent to a vertex \sout{$V$} of $C$ {\bf where at least two other
% edges $e_1$ and $e_2$ emanate from as indicated in the picture (due
% to the fact that a tropical curves does not contain $2$-valent
% vertices).} {\bf (Alternativ koennte man $V$ ins Bild einzeichnen
% wenn man dem Vertex einen Namen geben will - allerdings hab ich die
% fig-datei nicht)}.
%  \sout{. As $V$ is at
% least $3$-valent, there are at least two other edges $e_1$ and $e_2$
% emanating from $V$ as indicated in the picture.} \sout{The picture
% shows an example.}
% \begin{center}
% \input{linear.pstex_t}
% \end{center}
% {\bf Ich finde ''reached'' im Kontext nicht wirklich definiert...wie
% waers mit etwas wie''let $a,b\in [n]$ be leaves in different
% connected components of $C\setminus\{k\}$ (there are at least two,
% due to the fact that $2$-valent vertices do not exist)''
% --- dann brauchen wir auch die Kanten $e_1$ und $e_2$ nicht mehr und koennten
% moeglicherweise auf das Bild verzichten.}
% 
%  Let $a$ be a marked end which can be reached from $V$ via $e_1$
% and $b$ a marked end which can be reached via $e_2$. 
Then
$T:=\{a,b\}$ is not contained in any of the sets $I_i$. Hence
$\lambda_{i,T}=0$ for all $i$. In particular $\sum_{i=1}^r \mu_i
\lambda_{i,T}=0$ and the equation \ref{eq-pos-rep} is a positive
representation. Therefore
\begin{displaymath}f_k(p)=f_k\left(\sum_{v_S\in
V_k}\big(\sum_{i=1}^r \mu_i \lambda_{i,S}\big) v_S\right)
=\sum_{v_S\in V_k}\big(\sum_{i=1}^r \mu_i \lambda_{i,S}\big)
f_k(v_S)\end{displaymath} and $f_k$ is linear on $\tau$.
\end{proof}

\begin{remark}\label{rem-order}
Let us explain in more detail how to compute $\divi(f_k)$ for a $d$-cycle $Z$ of
$\calM_{0,n}$. We require that $Z$ is supported on the cones of
$\calM_{0,n}$ corresponding to combinatorial types. As $f_k$ is
linear on each cone by lemma \ref{def-fk}, the locus of
non-differentiability of $f_k$ is contained in cones of codimension
$1$ of $Z$. Given a cone $\tau$ of codimension $1$ in $Z$, we need
to compute the weight $\omega(\tau)$ of the cone $\tau$ \begin{displaymath}\omega(\tau)=\sum_i f_k\left(\omega(\sigma_i)
u_{\sigma_i/\tau}\right)- f_k\left(\sum_i \omega(\sigma_i)
u_{\sigma_i/\tau}\right)\end{displaymath} where $\sigma_i$ denote
the top-dimensional neighboring cones of $\tau$, $\omega(\sigma_i)$
their weight and $u_{\sigma_i/\tau}$ their primitive vectors. In
our case the primitive vectors are given by the structure of a marked
fan: each primitive vector corresponds to a to a tropical curve
with only one bounded edge of length $1$. Each $u_{\sigma_i/\tau}$
is equal to $v_{I_i}$ for a subset $I_i\subset [n]$ (assume $k\notin
I_i$ for all $i$). (The marked leaves in $I_i$ are on one side of
the bounded edge.) As $f_k$ is defined on the set $V_k$, we need to
find a positive representation of the $u_{\sigma_i/\tau}$ with
respect to $V_k$ (see \ref{rem-basisdarstellung} and \ref{rem-basisdarstellung2}). This can
be done using lemma \ref{basisrepresentation}. Also, we need to find
a positive representation of $\sum_i
\omega(\sigma_i)u_{\sigma_i/\tau}$. Given a representation
$u_{\sigma_i/\tau}=\sum_{T\in V_k} \lambda_{i,T} v_T$ (note that all
$\lambda_{i,T}$ are either $0$ or $1$ by the above), the following
equality holds
\begin{displaymath}
\sum_i \omega(\sigma_i) u_{\sigma_i/\tau}= \sum_i \omega(\sigma_i) \big(\sum_{T\in V_k} \lambda_{i,T} v_T\big) = \sum_{T\in V_k} \big(\sum_i \omega(\sigma_i) \lambda_{i,T}\big) v_T.
\end{displaymath}
Of course $\sum_i \omega(\sigma_i) \lambda_{i,T}\geq 0$ for all $T$, but the equation above is not necessarily a positive representation, since it is possible that none of the coefficients is zero. To make it a positive representation, we have to subtract $\sum_{T\in V_k} \min_{T\in V_k} (\sum_i  \omega(\sigma_i)\lambda_{i,T}) v_T$.
Then
\begin{align*}
&\sum_i f_k( \omega(\sigma_i) u_{\sigma_i/\tau})- f_k\big(\sum_i \omega(\sigma_i) u_{\sigma_i/\tau}\big)\\=& \sum_i \sum_{T\in V_k} \omega(\sigma_i) \lambda_{i,T} - \big( \sum_{T\in V_k} (\sum_i  \omega(\sigma_i)\lambda_{i,T}) - \sum_{T\in V_k} \min_{T\in V_k} (\sum_i  \omega(\sigma_i)\lambda_{i,T}) \big)\\=& \min_{T\in V_k} (\sum_i  \omega(\sigma_i) \lambda_{i,T}) \cdot \binom{n-1}{2},
\end{align*}
because $|V_k|= \binom{n-1}{2}$.
That is, in order to determine the weight of the cone $\tau$ in $\divi(f_k)$, we only have to determine the value $\min_{T\in V_k} (\sum_i  \omega(\sigma_i) \lambda_{i,T})$.
Recall that $\lambda_{i,T}=1$ if $T\subset I_i$ and $0$ else. Hence \begin{displaymath}\min_{T\in V_k} \big(\sum_i  \omega(\sigma_i) \lambda_{i,T}\big)= \min_{T\in V_k}\big(\sum_{i\;:\; T\subset I_i} \omega(\sigma_i)\big).\end{displaymath}
\end{remark}

\begin{proposition}\label{prop-basiccase}
Let $f_k$ be as in notation \ref{def-fk}. Then the
divisor of $f_k$ in $\calM_{0,n}$ is
\begin{displaymath}\divi(f_k)=\binom{n-1}{2}\Psi_k,\end{displaymath}
where $\divi(f_k)$ is defined in 3.4 of \cite{AR} and $\Psi_k$ is
defined in definition \ref{def-Psi}.
\end{proposition}
\begin{proof}
We may assume without restriction that $k=1$. By lemma
\ref{lem-fklinear}, the locus of non-differentiability of $f_1$ is
contained in the cones of codimension $1$ of $\calM_{0,n}$. A cone
$\tau$ of codimension $1$ in $\calM_{0,n}$ corresponds to the
combinatorial type of a tropical curve $C$ with one $4$-valent
vertex. Let $A_1$ denote the subset of $[n]$ consisting of the ends
which can be reached from this $4$-valent vertex via the adjacent
edge $e_1$, $A_2$ the subset which can be reached via $e_2$ and so
on. Without restriction $1\in A_1$ and $2\in A_2$. The three
neighboring cones $\sigma_1$, $\sigma_2$ and $\sigma_3$ of $\tau$
are given by the three possible resolutions of the $4$-valent
vertex in two $3$-valent vertices (i.e.\ the three possible ways to add a new edge which produces two $3$-valent vertices instead of the one $4$-valent).
\begin{center}
\includegraphics{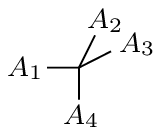}
\end{center}
\begin{center}
\includegraphics{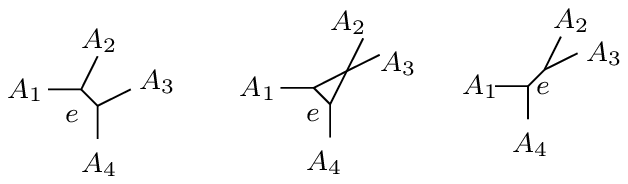}
\end{center}

The primitive vector $u_{\sigma_i/\tau}$ is given by the tropical curve where all edges except the new edge $e$ (which is introduced by resolving) are shrunk to length $0$. Hence the three primitive vectors are
\begin{align*}
v_{A_1\cup A_2}=v_{A_3 \cup A_4} \;,\;\;\; v_{A_1\cup A_3}=v_{A_2 \cup A_4} \;\;\mbox{ and}\;\;\;v_{A_1\cup A_4}=v_{A_2 \cup A_3}.
\end{align*}

By lemma \ref{basisrepresentation}, the positive representations of
these three primitive vectors with respect to $V_1$ (see \ref{rem-basisdarstellung}) are
\begin{displaymath} \sum_{T\subset A_3\cup A_4} v_T\;,\;\;\;\sum_{T\subset A_2\cup A_4} v_T\;\mbox{, respectively}\;\;\sum_{T\subset A_2\cup A_3} v_T,\end{displaymath}
(where in each case the sum goes over all $T\in V_k$).

By remark \ref{rem-order}, the weight of $\divi(f_1)$ along $\tau$
is equal to $\min_{T\in V_k} (|\{i\;:\;T\subset I_i\}| ) \cdot
\binom{n-1}{2}$, where $I_1=A_3 \cup A_4$, $I_2=A_2 \cup A_4$ and
$I_3=A_2 \cup A_3$ (since all weights of the neighbouring cones are
$1$). There are two cases to distinguish. Assume first that $\{1\}
\neq A_1$, that is, there exists a number $1\neq a_1\in A_1$. Then for
$T=\{a_1,2\}$ the number $|\{i\;:\;T\subset I_i\}|$ is zero. Hence the
weight of $\divi(f_1)$ at $\tau$ is zero, and $\tau$ is not part of
$\divi(f_1)$. Assume next that $A_1=\{1\}$. Then all vectors $v_T\in
V_1$ appear in the sum of the three primitive vectors. Let $a_2\in A_2$ and $a_3\in A_3$, then $T=\{a_2,a_3\}$
appears only once, that is, the number $|\{i\;:\;T\subset I_i\}|$ is
one. Hence the weight of $\divi(f_1)$ along $\tau$ is
$\binom{n-1}{2}$. In particular, $\divi(f_1)$ consists of the cones
given by a $4$-valent tropical curves such that leaf $1$ is
adjacent to the $4$-valent vertex, and each such cone has weight
$\binom{n-1}{2}$. It follows that $\divi(f_1)=\binom{n-1}{2}\cdot
\Psi_1$.
\end{proof}

  \section{Intersecting tropical Psi-classes}
In the sequel of this section, given integers $k_1,\ldots ,k_n\in\ZZ_{\geq 0}$ and a subset $I\subset [n]$, we denote $K(I)=\sum_{i\in I}k_i$. Using this notation, the main theorem of this section reads as follows:

\begin{theorem}\label{thm-psiintersect}
The intersection $\Psi_1^{k_1}\cdot\ldots\cdot \Psi_n^{k_n}$ is the
subfan of $\calM_{0,n}$ consisting of the closure of the cones of
dimension $n-3-K([n])$ corresponding to abstract tropical
curves $C$ such that for each vertex $V$ of $C$ we have $\val(V)=K(I_V)+3$, where $I_V$ denotes the set
$$I_V=\{i\in [n]\; : \text{ leaf }x_i \text{ is adjacent to } V \text{ and } k_i\geq 1\}\subset [n].$$
The weight of the facet $\sigma(C)$ containing the point $\varphi_n(C)$ equals
$$\omega(\sigma(C))=\frac{\prod_{V\in V(C)}K(I_V)!}{\prod_{i=1}^n k_i!}.$$
\end{theorem}

\begin{proof}
We prove by induction on $K([n])$, that is, we compute the weight $\omega(\tau)$ of a codimension one cell $\tau\subset \prod_{i=1}^n\psi_i^{k_i}$ in the intersection product $\psi_{1}\prod_{i=1}^n\psi_i^{k_i}$.
Let $\tau$ be a ridge of $\prod_{i=1}^n\psi_i^{k_i}$ and
let $C$ be a curve parameterized by $\tau$. As $\tau$ is of codimension one, there is by construction exactly one vertex $V$ of $C$ of valence one higher than expected, i.e.\ $\val(V)=K(I_V)+4$ --- apart from the leaves $x_i$ with $i\in I_V=\{i_1,\ldots ,i_l\}$ there are exactly $K(I_V)+4-|I_V|$ edges $a_1,\ldots a_t$ adjacent to $V$ as indicated in the picture below. 

\begin{center}
\includegraphics{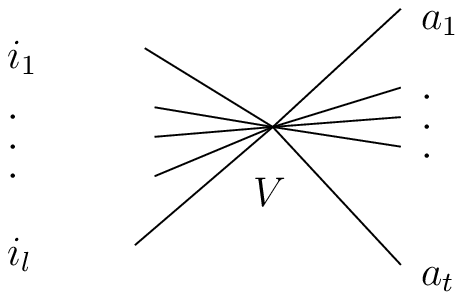}
\end{center}

%\begin{center}
%\input{chapter1/Pics/m052.pstex_t}
%\end{center}

% We let $A_k$ denote the
%set of all labels of leaves contained in the connected component of $C\setminus\{V\}$ containing the edge $a_k$.
Using remark \ref{rem-order}, we compute 
$$\omega(\tau)=\min_{v_T\in V_k}\left(\sum_{i:\;T\subset A_i} \omega(\sigma_i)\right)$$ % \cdot\binom{n-1}{2}$$
where $\sigma_i$ denote the facets containing $\tau$ and $v_{\sigma_i/\tau}=v_{A_i}$ denote their primitive vectors.

Facets $\sigma$ containing the ridge $\tau$ parameterize curves $C'$ such that $C$ is obtained from $C'$ by collapsing an edge $E$ with 
vertices $V_1$ and $V_2$ to the vertex $V$. By assumption, we know that the weight of such a facet equals
   $$\omega(\sigma_i)=\frac{\prod K(I_{V'})!}{\prod k_i!}\cdot K(I_{V_1})! K(I_{V_2})!=: W\cdot K(I_{V_1})! K(I_{V_2})!,$$
where the first product goes over all vertices $V'$ of $C'$ different from $V_1$ and $V_2$. The normal 
vector of such a cone $\sigma_i$ is by definition given by $v_{\sigma_i/\tau}=v_{A_i}$, where $A_i$ is the set
of labels on one of connecting component of $C'\setminus\{E\}$. Assume without loss of generality that $x_1$ is in the connected component of $C'\setminus\{E\}$ containing $V_1$ and let $A_i\subset [n]$ denote the subsets of labels of the other connected component.

\begin{center}
\includegraphics{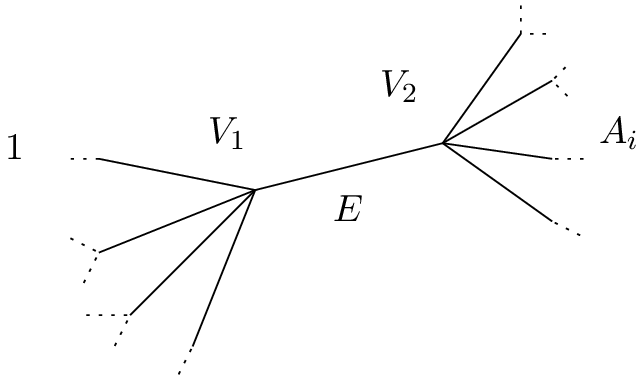}
\end{center}

Assume first that $x_1$ is not adjacent to $V$. Then there exist a leaf $x_s\neq x_1$ in the connected component of $C\setminus\{V\}$ containing $x_1$, hence the label $s$ is not contained in any of the sets $A_i$. 
For any $T$ containing $s$ we thus have $\sum_{i:\;T\subset A_i} \omega(\sigma_i)=0$ and hence the minimum over all $T$ is $0$, too. 
Consequently, 
$\omega(\tau)=0$.

Assume that $x_1$ is adjacent to $V$, let $s_1, s_2\in [n]\setminus\{1\}$ be two labels. It suffices to prove that 
$$
\sum_{T\subset A_i} \omega(\sigma_i) 
 =  \frac{\prod K(I_{V'})!}{\prod k_i!}\cdot \frac{(K(I_V)+1)!}{k_1+1}\\
 =  W \cdot \frac{(K(I_V)+1)!}{k_1+1}
$$
where the first product goes over all vertices $V'$ of $C$ different from $V$.

Note that a cell $\sigma$ satisfying the conditions above corresponds to a 
partition $J_{V_1}\cup J_{V_2}=I_V\setminus\{1,s_1,s_2\}=:M$ and a distribution 
of the labels $a_i$ among the leaves not labelled by some $x_i$, $i\in M\cup\{1,s_1,s_2\}$.
       %$i\in I_V\cup \{1,s_1,s_2\}$.
The total number of labels $a_i$ equals $K(I_V)+4-|M|$ % |I_{V}\cup \{1,s_1,s_2\}|$.
and there are exactly $K(J_{V_1})+k_1+1-|J_{V_1}|$ non-labelled edges at vertex $V_1$.
 %\item The 'weight' of each curve equals $(K(I_{V_1})+k_1)!\cdot (K(I_{V_2})+k_{s_1}+k_{s_2})!$.
Hence we get
\begin{eqnarray*}
\omega(\tau) & = & 
W\;\cdot \; \sum_{J_{V_1}\cup J_{V_2}=M}%I_V\setminus\{1,s_1,s_2\}} 
            (K(J_{V_1})+k_1)!(K(J_{V_2})+k_{s_1}+k_{s_2})!
            \binom{K(I_V)+1-|M|}{K(J_{V_1})+k_1+1-|J_{V_1}|}\\
\end{eqnarray*} 
and we want to see that this is equal to

\begin{displaymath}
 W \cdot \frac{(K(I_V)+1)!}{k_1+1}.
\end{displaymath}

The equality follows after proving the following identity
\begin{equation}\label{eq11}\sum_{I\subset M} 
         (K(I)+k_1)!(K-(K(I)+k_1))!
         \binom{K+1-m}{K(I)+k_1+1-|I|}=\frac{(K+1)!}{k_1+1}\end{equation}
setting $I=J_{V_1}$, $m= |M|$ and $K=K(I_V)$.

To see that this identity holds, multiply the left hand side by $\frac{k_1+1}{(K+1-m)!}$ and simplify:
\begin{align*}
&(k_1+1)\cdot\left(\sum_{I\subset M} 
         \frac{(K(I)+k_1)!\cdot(K-(K(I)+k_1))!}{(K(I)+k_1+1-|I|)!\cdot (K-K(I)-k_1-m+|I|)! }\right)\\
=& (k_1+1)\cdot\left( \frac{(k_1)! \cdot(K-k_1)!}{(k_1+1)!\cdot (K-k_1-m)!} +\sum_{\emptyset \neq I\subset M} 
         (K(I)+k_1)^{[|I|-1]} (K-K(I)-k_1)^{[m-|I|]}\right)\\
=& (K-k_1)^{[m]}+ \sum_{\emptyset \neq I\subset M} 
      (k_1+1)   (K(I)+k_1)^{[|I|-1]} (K-K(I)-k_1)^{[m-|I|]}
\end{align*}
where we use the falling power notation \begin{displaymath}x^{[p]}=x\cdot (x-1)\cdot \ldots\cdot (x-p+1).\end{displaymath}
If we multiply the right hand side of equation (\ref{eq11}) by $\frac{k_1+1}{(K+1-m)!}$ we get
$(K+1)^{[m]}$. Thus the identity follows from lemma \ref{applem2}.

\end{proof}
\begin{corollary}
If the intersection $\Psi_1^{k_1}\cdot\ldots\cdot \Psi_n^{k_n}$ is $0$-dimensional, i.e. $K([n])=n-3$,
then the (stable) intersection  $\Psi_1^{k_1}\cdot\ldots\cdot \Psi_n^{k_n}$ is just the origin $\{0\}$, with
weight $$\omega(\{0\})=\frac{(n-3)!}{k_1!\dots
k_n!}=\binom{n-3}{k_1,\dots ,k_n}.$$

\end{corollary}

\begin{remark}
For $n\geq 3$, let $\overline{M}_{0,n}$ denote the space of
$n$-pointed stable rational curves, that is, tuples $(C,p_1,\dots
,p_n)$ consisting of
 a connected algebraic curve $C$ of arithmetic genus $0$ with simple nodes as only
singularities and a collection $p_1,\dots ,p_n$ of distinct smooth
points on $C$ such that the number of automorphisms of $C$ with the
property that the points $p_i$ are fixed is finite.

 Define the line bundle $\mathcal{L}_i$ on
$\overline{M}_{0,n}$ to be the unique line bundle whose fiber over
each pointed stable curve $(C,p_1,\dots ,p_n)$ is the cotangent
space of $C$ at $p_i$ and let $\Psi_i\in A^1(\overline{M}_{0,n})$
denote its first Chern class. For $1\leq i\leq n$, let
$k_i\in\ZZ_{\geq 0}$ such that $\sum_{i=1}^n
k_i=\dim(\overline{M}_{0,n})=n-3$. Then the following equation holds
(see e.g. \cite{HM}, section 2.D)
\begin{equation*}
\int_{\overline{M}_{0,n}}\Psi_1^{k_1}\Psi_2^{k_2}\dots\Psi_n^{k_n}
=\frac{(n-3)!}{\prod_{i=1}^n k_i!}.
\end{equation*}
Hence the $0$-dimensional intersection products of Psi-classes
on the moduli space of $n$-marked rational algebraic curves
coincides with its tropical counterpart.
\end{remark}

\begin{example}[Psi-classes on $\mathcal{M}_{0,5}$]
If we intersect two $\Psi$-classes, the intersection is $0$-dimensional. Hence there are (up to symmetry) only two different intersection products to compute:
$\Psi_1^2$ and $\Psi_1\cdot \Psi_2$.
Let us compute both.
Let us start with $\divi(f_1)\cdot \Psi_1$. By lemma \ref{lem-fklinear} we know that we only have to check the cones of codimension $1$ in $\Psi_1$ --- that is, the cone $\{0\}$.
The neighbors of $\{0\}$ in $\Psi_1$ --- that is, in this case, the top-dimensional cones of $\Psi_1$ --- correspond to tropical curves with $1$ at a $4$-valent vertex:
\begin{center}
\includegraphics{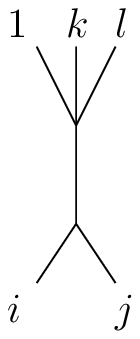}
\end{center}
(Here we assume that $\{i,j,k,l\}=\{2,3,4,5\}$.)
There are $\binom{4}{2}=6$ of these cones. Each such cone is generated by the primitive vector $v_{\{i,j\}}$. The sum over all primitive vectors is $0$.
Hence the weight of $\{0\}$ is given by
\begin{displaymath}
\sum_{i,j\in\{2,3,4,5\}, i\neq j}f_1(v_{\{i,j\}})- f_1(0)= 6.
\end{displaymath}
(By definition, all cones of $\Psi_1$ are of weight $1$.) Thus the
weight of $\{0\}$ in $\divi(f_1)\cdot \Psi_1$ is $6$, and using
proposition \ref{prop-basiccase}, the weight of $\{0\}$ in
$\Psi_1\cdot \Psi_1$ is $1$.

 Now let us compute the weight of $\{0\}$ in $\divi(f_1)\cdot
\Psi_2$. Three of the neighbors of $\{0\}$ correspond to a curve as
on the right, the other three to a curve as on the left:
\begin{center}
\includegraphics{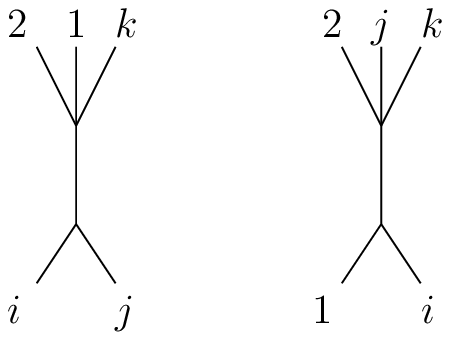}
\end{center}
(We assume $\{i,j,k\}=\{3,4,5\}$.)
The primitive vectors of the first type are $v_{i,j}$ and they are
given as a positive combination with respect to $V_1$. A positive
combination for the primitive vectors of the second type is
$v_{\{2,j\}}+v_{\{2,k\}}+v_{\{j,k\}}$. Their sum is of course again $0$. Hence
the weight of $\{0\}$ is
\begin{align*}
&f_1(v_{\{3,4\}})+f_1(v_{\{3,5\}})+f_1(v_{\{4,5\}})+f_1(v_{\{2,3\}}+v_{\{2,4\}}+v_{\{3,4\}})\\&+
f_1(v_{\{2,3\}}+v_{\{2,5\}}+v_{\{3,5\}})+f_1(v_{\{2,4\}}+v_{\{2,5\}}+v_{\{4,5\}}) -
f_1(0)\\=&1+1+1+3+3+3=12.
\end{align*}
Thus --- using proposition \ref{prop-basiccase} again --- the weight of $\{0\}$ in $\Psi_1\cdot \Psi_2$ is $2$.
\end{example}

\begin{example}[Psi-classes on $\mathcal{M}_{0,6}$]
Let us compute $\Psi_1^2$, respectively $\divi(f_1)\cdot \Psi_1$, on $\mathcal{M}_{0,6}$. To do so, we have to compute the weight of a cone of codimension $1$ in $\Psi_1$. Such a cone corresponds to a tropical curve with either another $4$-valent vertex (as on the left) or with a $5$-valent vertex, to which $1$ is adjacent (as on the right):
\begin{center}
\includegraphics{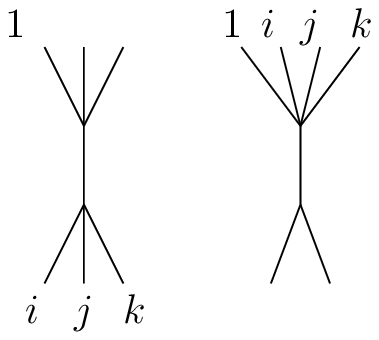}
\end{center}
A cone corresponding to the curve on the left has $3$ neighbors corresponding to the $3$ possible resolutions of the lower vertex. These three cones are generated by the primitive vectors $v_{\{i,j\}}$, $v_{\{i,k\}}$ and $v_{\{j,k\}}$. Thus the weight of such a cone is
\begin{displaymath}
f_1(v_{\{i,j\}})+f_1(v_{\{i,k\}})+f_1(v_{\{j,k\}})-f_1(v_{\{i,j\}}+v_{\{i,k\}}+v_{\{j,k\}})=0,
\end{displaymath}
and it does not belong to $\Psi_1^2$. A cone corresponding to
the curve on the right (where we assume now $\{i,j,k\}=\{2,3,4\}$
for simplicity) has $6$ neighbors in $\Psi_1$, $3$ as on the right
and $3$ as on the left (below the curve, the corresponding normal
vector for the cones are shown):
\begin{center}
\includegraphics{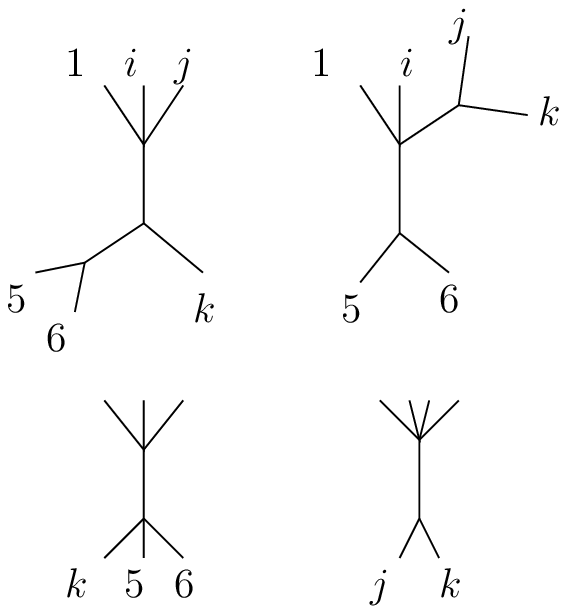}
\end{center}
Hence, the weight of this cone is
\begin{align*}
&f_1(v_{\{2,5\}}+v_{\{2,6\}}+v_{\{5,6\}})+ f_1(v_{\{3,5\}}+v_{\{3,6\}}+v_{\{5,6\}})+f_1(v_{\{4,5\}}+v_{\{4,6\}}+v_{\{5,6\}})\\&+f_1(v_{\{2,3\}})+f_1(v_{\{2,4\}})+f_1(v_{\{3,4\}})
\\&- f_1( v_{\{2,5\}}+v_{\{2,6\}}+v_{\{5,6\}}+v_{\{3,5\}}+v_{\{3,6\}}+v_{\{5,6\}}+v_{\{4,5\}}+v_{\{4,6\}}\\&\qquad\qquad\qquad+v_{\{5,6\}}+v_{\{2,3\}}+v_{\{2,4\}}+v_{\{3,4\}})
\\=&3+3+3+1+1+1-f_1(2 v_{\{5,6\}})=
12-2 =10,
\end{align*}
because the sum over all primitive vectors is not automatically written
as a positive combination with respect to $V_1$ and we have to
subtract $\sum_{i,j\in [6]\setminus\{1\},i\neq j }v_{\{i,j\}}$ to make a positive
combination. By proposition \ref{prop-basiccase}, the weight of each
such cone in $\Psi_1^2$ is one.

\end{example}

% \begin{example}
% If only one $\Psi$-class is involved $\Leftarrow$ 1\\
%
% If $k_1=n-4$, $k_2=1$, then after reordering $\Psi_2(\Psi_1^{n-4})$
% , all weights are trivial
% \end{example}

\begin{appendix}
 
\section{Equations needed for the proof the main theorem, \ref{thm-psiintersect}}

\normalfont
Let $k_1,\ldots,k_m$ and $K$ be integers. For a subset $I\subset [m]:=\{1,\ldots,m\}$ we use the notation $K(I)=\sum_{i\in I}k_i$ and the falling power notation \begin{displaymath}x^{[p]}=x\cdot (x-1)\cdot \ldots\cdot (x-p+1).\end{displaymath}

\begin{lemma}\label{applem1}
The following equation is satisfied:
\begin{equation}\label{fpi1}
\sum_{\emptyset \neq I\subset [m]} K(I)^{[|I|-1]}(K-K(I))^{[m-|I|]}=m\cdot K^{[m-1]}.
\end{equation}
\end{lemma}
\begin{proof}

The proof is an induction. We assume that equation (\ref{fpi1}) holds for any $p<m$ and prove that it holds for $m$. The induction beginning holds trivially.

We use the well-known multinomial identity for falling powers:
\begin{equation}\label{fpi2}
(x_1+\ldots+x_n)^{[p]}= \sum_{a_1+\ldots+a_r=p}\frac{p!}{a_1!\cdot \ldots \cdot a_r!}x_1^{[a_1]}\cdot \ldots\cdot x_r^{[a_r]}
\end{equation}
where the $a_i$ are nonnegative integers.
Set $x_i=k_i$ for all $i=1,\ldots,m$ and $x_0=K-\sum_{i\in [m]} k_i$.
Express each side of equation (\ref{fpi1}) as a linear combination of monomials $x_0^{[a_0]}\cdot \ldots\cdot x_r^{[a_r]}$ with $a_0+\ldots+a_r=m-1$ and compare the coefficients.
The coefficient of such a monomial on the right hand side of (\ref{fpi1}) equals $m\cdot\frac{(m-1)!}{a_0!\ldots a_r!}=\frac{m!}{a_0!\ldots a_r!}$. The coefficient on the left hand side equals
\begin{displaymath}
\sum_{I}\frac{(|I|-1)!(m-|I|)!}{a_0!\ldots a_r!}
\end{displaymath}
where the sum is over all non-empty subsets of $[m]$ satisfying $\sum_{i\in I}a_i= |I|-1$.
To prove (\ref{fpi1}) we thus have to show the following identity for any tuple $(a_1,\ldots,a_m)$ of nonnegative integers satisfying $a_1+\ldots+a_m<m$:
\begin{equation}\label{fpi3}
\sum_I (|I|-1)!(m-|I|)! = m!
\end{equation}
where again the sum goes over all $I\subset [m]$ satisfying $\sum_{i\in I}a_i= |I|-1$.
We use induction on $\sum_{i\in [m]}a_i$ to prove equation (\ref{fpi3}). If all $a_i=0$ then the only subsets $I\subset [m]$ satisfying $\sum_{i\in I}a_i= |I|-1$ are one-element subsets. There are $m$ of those and they all contribute a summand of $(m-1)!$ to the left hand side, so we have $m!$ altogether which equals the right hand side.
Now assume that $a_1,\ldots,a_p$ are positive and $a_{p+1}=\ldots=a_m=0$ for some $p$. Then $p<m$ since $a_1+\ldots+a_m<m$. We can write any subset $I\subset [m]$ as a disjoint union $I=J\cup I'$ where $J\subset [p]$ and $I'\subset \{p+1,\ldots,m\}$.
The sum $\sum_{i\in I}a_i$ depends only on $J$, i.e.\ it is equal to $\sum_{i\in J}a_i$. As before we denote $ \sum_{i\in J}a_i$ by $K(J)$. If we fix a set $J$ we can produce several possible subsets $I$ satisfying $K(J)=\sum_{i\in J}a_i=\sum_{i\in I}a_i= |I|-1$ by just adding $K(J)-|J|+1$ elements of $\{p+1,\ldots,m\}$. Therefore we can write the left hand side of equation (\ref{fpi3}) as
\begin{align*}
&\sum_{J\subset [p]} \binom{m-p}{K(J)-|J|+1} K(J)! (m-1-K(J))!\\
=& (m-p) (m-1)!+\sum_{\emptyset \neq J\subset [p]}\frac{(m-p)! K(J)! (m-1-K(J))! }{(K(J)-|J|+1)!(m-p-K(J)+|J|-1)!}\\
=& (m-p) (m-1)!+(m-p)!\sum_{\emptyset \neq J\subset [p]} K(J)^{[|J|-1]} (m-1-K(J))^{[m-|J|]}.
\end{align*}
Subtract the contribution of the empty set from the right hand side of equation (\ref{fpi3}) and divide by $(m-p)!$. Then we get
\begin{align*}
\frac{m!- (m-p) (m-1)!}{(m-p)!}= \frac{(m-1)!}{(m-p)!}(m-m+p)= p\cdot (m-1)^{[p-1]}. 
\end{align*}
Hence (\ref{fpi3}) follows if
\begin{displaymath}
\sum_{\emptyset \neq J\subset [p]} K(J)^{[|J|-1]} (m-1-K(J))^{[m-|J|]}=p\cdot (m-1)^{[p-1]}
\end{displaymath}
which holds by the induction assumption on (\ref{fpi1}).
\end{proof}

\begin{lemma}\label{applem2}
Let $M=\{2,\ldots , m+1\}$ and let $K=\sum_{i\in [m+1]}k_i$. 
Then the following equation holds:
\begin{equation}\label{eq1}
(K-k_1)^{[m]}+\sum_{\emptyset\neq I\subset M} (k_1+1)\cdot (K(I)+k_1)^{[|I|-1]}(K-K(I)-k_1)^{[m-|I|]}=(K+1)^{[m]}.
\end{equation}

\end{lemma}
\begin{proof}
The proof is a double induction on $k_1$ and $m$. We show that the equation is true for $k_1=0$ and for all $m$. Next, we assume that it is true for all $k_1-1$ and any $m$ and for any $k_i$ and $m-1$ and show that it is true for $k_1$ and $m$.
For $k_1=0$, the equation reads
\begin{displaymath}
K^{[m]}+\sum_{\emptyset\neq I\subset M}  (K(I))^{[|I|-1]}(K-K(I))^{[m-|I|]}=(K+1)^{[m]}.
\end{displaymath}
If we subtract $K^{[m]}$ from the right hand side, we get \begin{align*}&(K+1)^{[m]}-K^{[m]}\\=& \big((K+1)-(K-m+1)\big) (K\cdot \ldots\cdot (K-m+2))\\=&m\cdot K^{[m-1]}.\end{align*}
Thus the equation for $k_1=0$ follows from lemma \ref{applem1} after relabeling the index set $M$.
Now we assume that the equation is true for $k_1-1$. 
Remember that $K$ is defined as $K=k_1+\sum_{i\in M}k_i$, so if we replace $k_1$ by $k_1-1$ then we also have to replace $K$ by $K-1$.
Then the equation reads
\begin{equation}\label{eq2}
(K-k_1)^{[m]}+\sum_{\emptyset\neq I\subset M} k_1\cdot (K(I)+k_1-1)^{[|I|-1]}(K-K(I)-k_1)^{[m-|I|]}=K^{[m]}.
\end{equation}
We subtract the left hand side of equation (\ref{eq2}) from the left hand side of equation (\ref{eq1}) and get
\begin{align*}
&\sum_{\emptyset\neq I\subset M} (k_1+1)\cdot (K(I)+k_1)^{[|I|-1]}(K-K(I)-k_1)^{[m-|I|]}
\\& - \sum_{\emptyset\neq I\subset M} k_1\cdot (K(I)+k_1-1)^{[|I|-1]}(K-K(I)-k_1)^{[m-|I|]}
\\=& \sum_{I=\{j\},j\in M}(K-k_j-k_1)^{[m-1]}\\&
+\sum_{|I|\geq 2} \big((k_1+1)(K(I)+k_1)-k_1(K(I)+k_1-|I|+1)\big)\\&\;\;\;\cdot(K(I)+k_1-1)^{[|I|-2]}(K-K(I)-k_1)^{[m-|I|]}
\end{align*}
which can be simplified to
\begin{multline}\label{eq4}
  \sum_{I=\{j\},j\in M}(K-k_j-k_1)^{[m-1]}\\
+\sum_{|I|\geq 2} (K(I)+k_1\cdot |I|) (K(I)+k_1-1)^{[|I|-2]}(K-K(I)-k_1)^{[m-|I|]}.
\end{multline}
Now we subtract the right hand side of equation (\ref{eq2}) from the right hand side of equation (\ref{eq1}) and get
\begin{equation}\label{eq3}
(K+1)^{[m]}-K^{[m]}= m\cdot K^{[m-1]}=\sum_{j\in M} K^{[m-1]}.
\end{equation}
We want to apply equation (\ref{eq1}) for $m-1$ and any $k_i$ differently for each of the $m$ summands above. To do so, we need to interpret $K-1$ as a sum of $m$ numbers. We choose the first summand (i.e.\ the analogue of $k_1$) to be $k_1+k_j-1$ and the other summands to be the $k_i$ (except $k_j$).
If we replace $k_1$ by $k_1+k_j-1$, $K$ by $K-1$ and $m$ by $m-1$ in equation (\ref{eq1}) it reads

\begin{align*}
&(K-k_1-k_j)^{[m-1]}\\&+\sum_{\emptyset\neq J\subset M\setminus\{j\}} (k_1+k_j)\cdot (K(J)+k_1+k_j-1)^{[|J|-1]}(K-K(J)-k_1-k_j)^{[m-|J|-1]}=K^{[m-1]}.
\end{align*}
Thus equation (\ref{eq3}) equals
\begin{multline}\label{eq5}
\sum_{j\in M}(K-k_1-k_j)^{[m-1]}\\+\sum_{\emptyset\neq J\subset M\setminus\{j\}} (k_1+k_j)\cdot (K(J)+k_1+k_j-1)^{[|J|-1]}(K-K(J)-k_1-k_j)^{[m-|J|-1]}.
\end{multline}
It remains to show that the expression (\ref{eq4}) equals expression (\ref{eq5}). To see this, note that every $I\subset M$ yields a possible $J$ for every $j\in I$ by just deleting $j$. Thus $|I|-1=|J|$. Any $I$ contributes a factor of $(k_1+k_j)$ in the summand for $j$ in (\ref{eq5}). Thus it contributes $K(I)+k_1\cdot |I|$ in total, which equals the contribution in (\ref{eq4}).
\end{proof}

\end{appendix}

\end {document}